\documentclass[12pt]{article}
\usepackage{amsmath}
\usepackage{amssymb}
\usepackage{amsfonts}
\usepackage{myart}

\input tcilatex
\begin{document}

\title{Monistic conception of geometry}
\author{Yuri A.Rylov}
\date{Institute for Problems in Mechanics, Russian Academy of Sciences,\\
101-1, Vernadskii Ave., Moscow, 119526, Russia.\\
e-mail: rylov@ipmnet.ru\\
Web site: {$http://rsfq1.physics.sunysb.edu/\symbol{126}rylov/yrylov.htm$}\\
or mirror Web site: {$http://gasdyn-ipm.ipmnet.ru/\symbol{126}%
rylov/yrylov.htm$}}
\maketitle

\begin{abstract}
One considers the monistic conception of a geometry, where there is only one
fundamental quantity (world function). All other geometrical quantities a
derivative quantities (functions of the world function). The monisitc
conception of a geometry is compared with pluralistic conceptions of a
geometry, where there are several independent fundamental geometrical
quantities. A generalization of a pluralistic conception of the proper
Euclidean geometry appears to be inconsistent, if the generalized geometry
is inhomogeneous. In particular, the Riemannian geometry appears to be
inconsistent, in general, if it is obtained as a generalization of the
pluralistic conception of the Euclidean geometry.
\end{abstract}

\section{Introduction}

The term "monism" (monistic) in its application to geometry means, that the
conception of a geometry uses only one fundamental quantity (distance) and
all other geometrical quantities and concepts are derivatives of the
fundamental quantity. The pluralistic conception is a conception, containing
several fundamental quantities, which may be independent. In general, a
monistic conception can be considered as a pluralistic conception. It is
sufficient to declare, that some derivative quantities of monistic
conception are fundamental quantities of a pluralistic conception. However,
not any pluralistic conception can be considered as a monistic conception.
It is necessary for such a consideration, that there are proper connections
between fundamental quantities of the pluralistic conception.

If one is going to modify or to generalize a conception, it is desirable to
present the conception in the form of a monistic conception. In this case it
is sufficient to modify properties of the only fundamental quantity. Other
(derivative) quantities will be modified automatically, because other
quantities are functions of the unique fundamental quantity. If the modified
conception is pluralistic, the modification becomes to be difficult, because
in the pluralistic conception there may be connections between the
fundamental quantities. In this case one may not modify independently
properties of different fundamental quantities, which look as independent.

Usually the monistic conception is a more developed conception, than the
foregoing pluralistic conception, containing several fundamental quantities.
For instance, the Christianity, which contains only one god, is a more
developed conception, than the heathendom, which contains many gods. Usually
a monistic conception (or a less pluralistic conception) is a result of
development of a pluralistic conception, which is followed by a reduction of
fundamental quantities. For instance, theory of thermal phenomena developed
from the axiomatic thermodynamics to the statistical physics. This
development was followed by transformation of fundamental thermodynamic
quantities into derivative quantities of the statistical physics. Appearance
of intermediate derivative concepts of the monistic conception makes the
monistic conception to be a more complicated conception, than the preceding
pluralistic conception. This complexity is conditioned by a necessity of
deduction of transformed quantities, which were fundamental quantities in
the pluralistic conception. However, the more complicated monistic
conception admits one to describe and explain such physical phenomena, which
cannot be described in the framework of the pluralistic theory. For
instance, the thermal fluctuations are described by the statistical physics,
but they cannot be described in the framework of the axiomatic
thermodynamics.

Geometry is a science on disposition and shape of geometrical objects in the
space or in the event space (space-time). Any geometry is given on a set $%
\Omega $ of points (events). Any geometric object $\mathcal{O}$ is a set $%
\Omega ^{\prime }$ of points $P\in \Omega ^{\prime }$ with ($\Omega ^{\prime
}\subset \Omega $). The shape of the geometric object $\mathcal{O}$ is
described, if the distance $\rho \left( P_{1},P_{2}\right) ,$ $\forall
P_{1},P_{2}\in {\Omega }^{\prime }$ is given. The shape and mutual
dispositions of all geometrical objects are given completely, if the
distance $\rho \left( P_{1},P_{2}\right) ,$ $\forall P_{1},P_{2}\in \Omega $
is given. In the usual space the distance $\rho $ is a real nonnegative
quantity. In the event space the distance $\rho $ is either a real
nonnegative quantity, or an imaginary quantity. In this case it is more
convenient to use the quantity $\sigma =\frac{1}{2}\rho ^{2}$, which is
always real. The quantity $\sigma $ is called the world function. This
quantity was introduced by J.L. Synge \cite{S60} for description of the
Riemannian geometry.

Thus, the geometry of the event space (space-time geometry) is described
completely by one real function $\sigma $.

\begin{equation}
\sigma :\quad \Omega \times \Omega \mathcal{\rightarrow }\mathbb{R},\quad
\sigma \left( P,Q\right) =\sigma \left( Q,P\right) ,\quad \sigma \left(
P,P\right) =0,\quad \forall P,Q\in \mathcal{\Omega }  \label{a1.1}
\end{equation}%
As far as a motion of physical bodies in the space is described as a
geometrical object in the space-time, a construction of geometrical objects
and investigation of their properties is very important in such applied
sciences as physics and mechanics. Any geometrical object $\mathcal{O}$ is
described by enumeration of points of the object $\mathcal{O}$. However,
such a description by means of a direct enumeration of points, belonging to $%
\mathcal{O}$, contains too much information. One tends to simplify this
procedure, constructing geometric objects from standard blocks. These blocks
contain many points of the space-time, and construction of the geometric
object $\mathcal{O}$ is reduced to a finite or countable number of
operations with blocks.

Procedure of the geometric objects construction from blocks has been
developed in the proper Euclidean geometry $\mathcal{G}_{\mathrm{E}}$. It
can be used for description of usual space (but not for space-time), because
the world function $\sigma _{\mathrm{E}}$ of the geometry $\mathcal{G}_{%
\mathrm{E}}$ is nonnegative. The nonnegativity of $\sigma _{\mathrm{E}}$ is
used essentially at construction of geometrical objects in $\mathcal{G}_{%
\mathrm{E}}$. A direct enumeration of points of the geometrical object may
be replaced by some relations, which are fulfilled for points $P\in \mathcal{%
O}$ and only for them. However, such a description needs a use of some
numerical functions, given on $\Omega $, for instance, world function $%
\sigma $.

In reality the proper Euclidean geometry $\mathcal{G}_{\mathrm{E}}$ is
constructed usually from three kinds of blocks (point, segment of straight,
angle) without mention of the fact, that these blocks may be described in
terms of the world function $\sigma _{\mathrm{E}}$ and only in terms of $%
\sigma _{\mathrm{E}}$ (See details in \cite{R2007}). Properties of the three
blocks are postulated in the form of axioms. Using the rules of construction
of a geometric object from these blocks, one can formulate the proper
Euclidean geometry in the form of a logical construction. Basic statements
of the geometry have the form of geometrical axioms, formulated in terms of
basic geometric quantities (point, segment of straight and angle). Usually
the distance $\rho _{\mathrm{E}}$ (or world function $\sigma _{\mathrm{E}}$)$%
\ $is considered to be some derivative concept, which may be not mentioned
at all.

The fact, that the proper Euclidean geometry $\mathcal{G}_{\mathrm{E}}$ may
be formulated in terms of the world function (distance) is well known.
Nevertheless the distance is not used usually as a fundamental quantity of
the proper Euclidean geometry $\mathcal{G}_{\mathrm{E}}$. This quantity
admits one to \textit{formulate the geometry }$G_{\mathrm{E}}$\textit{\ as a
monistic conception, based on the only fundamental concept} (world function $%
\sigma _{\mathrm{E}}$). Of course, the world function $\sigma _{\mathrm{E}}$
of the geometry $\mathcal{G}_{\mathrm{E}}$ is to satisfy some conditions,
formulated in terms of the world function. The necessary and sufficient
conditions, that the geometry $\mathcal{G}$, described by the world function 
$\sigma $, is $n$-dimensional proper Euclidean geometry, have the form \cite%
{R2002}

\noindent I. 
\begin{equation}
\exists \mathcal{P}^{n}\subset \Omega ,\qquad F_{n}(\mathcal{P}^{n})\neq
0,\qquad F_{n+1}(\Omega ^{n+2})=0,  \label{a1.2}
\end{equation}%
where $\mathcal{P}^{n}=\left\{ P_{0},P_{1},..P_{n}\right\} $, the quantity $%
F_{n}(\mathcal{P}^{n})$ is the Gramm's determinant%
\begin{equation}
F_{n}(\mathcal{P}^{n})=\det \left\vert \left\vert \left( \mathbf{P}_{0}%
\mathbf{P}_{i}.\mathbf{P}_{0}\mathbf{P}_{k}\right) \right\vert \right\vert
,\quad i,k=1,2,...n  \label{a1.3}
\end{equation}%
and $\left( \mathbf{P}_{0}\mathbf{P}_{i}.\mathbf{P}_{0}\mathbf{P}_{k}\right) 
$ is the scalar product of two vectors $\mathbf{P}_{0}\mathbf{P}_{i}$ and $%
\mathbf{P}_{0}\mathbf{P}_{k}$, which is determined in terms of the world
function $\sigma $ by means of the relation%
\begin{equation}
\left( \mathbf{P}_{0}\mathbf{P}_{i}.\mathbf{P}_{0}\mathbf{P}_{k}\right)
=\sigma \left( P_{0},P_{k}\right) +\sigma \left( P_{0},P_{i}\right) -\sigma
\left( P_{i},P_{k}\right)  \label{a1.4}
\end{equation}

\noindent II. 
\begin{equation}
\sigma (P,Q)={\frac{1}{2}}\sum_{i,k=1}^{n}g^{ik}(\mathcal{P}%
^{n})[x_{i}\left( P\right) -x_{i}\left( Q\right) ][x_{k}\left( P\right)
-x_{k}\left( Q\right) ],\qquad \forall P,Q\in \Omega ,  \label{a1.5}
\end{equation}%
where the quantities $x_{i}\left( P\right) $, $x_{i}\left( Q\right) $ are
defined by the relations 
\begin{equation}
x_{i}\left( P\right) =\left( \mathbf{P}_{0}\mathbf{P}_{i}.\mathbf{P}_{0}%
\mathbf{P}\right) ,\qquad x_{i}\left( Q\right) =\left( \mathbf{P}_{0}\mathbf{%
P}_{i}.\mathbf{P}_{0}\mathbf{Q}\right) ,\qquad i=1,2,...n  \label{a1.6}
\end{equation}%
The contravariant components $g^{ik}(\mathcal{P}^{n}),$ $(i,k=1,2,\ldots n)$
of metric tensor are defined by its covariant components $g_{ik}(\mathcal{P}%
^{n}),$ $(i,k=1,2,\ldots n)$ by means of relations 
\begin{equation}
\sum_{k=1}^{n}g_{ik}(\mathcal{P}^{n})g^{kl}(\mathcal{P}^{n})=\delta
_{i}^{l},\qquad i,l=1,2,\ldots n  \label{a1.7}
\end{equation}%
where 
\begin{equation}
g_{ik}(\mathcal{P}^{n})=\left( \mathbf{P}_{0}\mathbf{P}_{i}.\mathbf{P}_{0}%
\mathbf{P}_{k}\right) ,\qquad i,k=1,2,\ldots n  \label{a1.8}
\end{equation}%
III.\quad The relations 
\begin{equation}
\left( \mathbf{P}_{0}\mathbf{P}_{i}.\mathbf{P}_{0}\mathbf{P}\right)
=x_{i},\qquad x_{i}\in \mathbb{R},\qquad i=1,2,\ldots n,  \label{a1.9}
\end{equation}%
considered to be equations for determination of $P\in \Omega $ as a function
of $x_{i},$ $i=1,2,...n$, have always one and only one solution.

The condition I determines the dimension $n$ of the geometry and $n$ basic
vectors $\mathbf{P}_{0}\mathbf{P}_{i}$, $i=1,2,...n$ on the set $\Omega $.
The condition II determines properties of the metric tensor $g_{ik}(\mathcal{%
P}^{n})$.

In such a representation the proper Euclidean geometry looks as a monistic
conception, which is described by the only fundamental quantity: world
function $\sigma $. All other geometrical quantities (concepts and
geometrical objects) are derivative in the sense, that all they can be
expressed in terms of points of the set $\Omega $ and in terms of world
functions between them.

Perspective of constructing such a monistic conception not only for the
proper Euclidean geometry seemed to be attractive. One introduced metric
space and constructed so called metric geometry. However, construction of
the metric geometry is possible only for the real metric (distance) $\rho
\geq 0$. Such a geometry cannot be constructed in the space-time. Besides,
for construction of straight lines (the shortest), one needs to impose on
the metric $\rho $ an additional condition (the triangle axiom) 
\begin{equation}
\rho \left( P,Q\right) +\rho \left( P,R\right) \geq \rho \left( R,Q\right)
,\quad \forall P,Q,R\in \Omega  \label{a1.10}
\end{equation}%
The triangle axiom (\ref{a1.10}) reflects our belief, that there are only
such geometries, where the straight is a one-dimensional point set.

Attempts \cite{M28,B53} were made to construct so called distance geometry,
which is free of the condition (\ref{a1.10}). Unfortunately, these attempts
failed in the sense, that the obtained distance geometry was not monistic
and completely metric. It contains also nonmetric procedure, which admits
one to construct straight lines by means of a continuous mapping of segment $%
[0,1]$ onto the set $\Omega $.

\section{Pluralistic approach to the Riemannian \newline
geometry as a reason of its inconsistency}

Usually one does not use the monistic conception of the proper Euclidean
geometry $\mathcal{G}_{\mathrm{E}}$, where there is the only fundamental
geometric quantity (world function). Instead, one uses pluralistic
conception of geometry, where there are several fundamental geometric
quantities (for instance, dimension, straight, angle,...), which are
considered as independent geometrical quantities. In the proper Euclidean
geometry one succeeds to agree properties of these different fundamental
quantities and to construct a consistent geometry. However, in the
non-Euclidean geometry one succeeds to agree properties of different
fundamental quantities not always. Having several fundamental concepts and
attributing to them some properties, one cannot be sure, that these
properties may be made compatible between themselves. As a result the
obtained geometry appears to be inconsistent. The problem of agreement
between properties of different fundamental quantities in a pluralistic
conception of geometry are formulated usually as a consistency of axioms of
the geometry. Any test of the geometric axioms consistency needs a lot of
labour, which can be carried out not always. As a result a pluralistic
geometry appears to be inconsistent. For instance, the Riemannian geometry,
which is constructed usually as a pluralistic conception, appears to be
inconsistent.

Using a monistic conception, when there is only one fundamental quantity,
one constructs a consistent geometry automatically. As an example, let us
consider construction of a straight line in a monistic geometry $\mathcal{G}$%
, described by the world function $\sigma $, given on the set $\Omega $ of
points $P$. We shall use term $\sigma $-space for the quantity $V=\left\{
\sigma ,\Omega \right\}$. The Euclidean space is a special case of the $%
\sigma $-space with $\sigma =\sigma _{\mathrm{E}}.$ In the proper Euclidean
geometry (Euclidean $\sigma $-space $V_{\mathrm{E}}=\left\{ \sigma _{\mathrm{%
E}},\Omega \right\} $), the segment $\mathcal{T}_{\left[ P_{0}P_{1}\right] }$
of straight line between the points $P_{0}$ and $P_{1}$ is defined by the
relation 
\begin{equation}
\mathcal{T}_{\left[ P_{0}P_{1}\right] }=\left\{ R|\sqrt{2\sigma \left(
P_{0},R\right) }+\sqrt{2\sigma \left( P_{1},R\right) }-\sqrt{2\sigma \left(
P_{0},P_{1}\right) }=0\right\}  \label{a1.12}
\end{equation}%
where the world function $\sigma $ coincides with the Euclidean world
function $\sigma _{\mathrm{E}}$. The same form of the straight segment $%
\mathcal{T}_{\left[ P_{0}P_{1}\right] }$ is to have in any other $\sigma $%
-space. In the $n$-dimensional $\sigma $-space the segment $\mathcal{T}_{%
\left[ P_{0}P_{1}\right] }$ is a $\left( n-1\right) $-dimensional surface,
in general. However, in the proper Euclidean geometry the $\left( n-1\right) 
$-dimensional surface degenerates into one-dimensional set. This fact is a
corollary of special properties (\ref{a1.2}) - (\ref{a1.9}) of the Euclidean
world function $\sigma _{\mathrm{E}}$, which generates fulfillment of the
triangle axiom (\ref{a1.10}). In the proper Riemannian geometry the triangle
axiom (\ref{a1.10}) is also valid. It is also connected with the form of the
Riemannian world function $\sigma _{\mathrm{R}}$, although in the
conventional (pluralistic) presentation of the Riemannian geometry, the
world function $\sigma _{\mathrm{R}}$ is a derivative quantity, defined by
the relation%
\begin{equation}
\sigma _{\mathrm{R}}\left( P_{0},P_{1}\right) =\frac{1}{2}\left(
\int_{P_{0}}^{P_{1}}\sqrt{g_{ik}\left( x\right) dx^{i}dx^{k}}\right) ^{2}
\label{a1.14}
\end{equation}%
where integral in (\ref{a1.14}) is taken along the geodesic, connecting
points $P_{0}$ and $P_{1}$.

World function $\sigma _{\mathrm{R}}$, defined by the relation (\ref{a1.14})
satisfies the equation \cite{S60}:%
\begin{equation}
\frac{\partial }{\partial x^{i}}\sigma _{\mathrm{R}}\left( x,x^{\prime
}\right) g^{ik}\left( x\right) \frac{\partial }{\partial x^{k}}\sigma _{%
\mathrm{R}}\left( x,x^{\prime }\right) =2\sigma _{\mathrm{R}}\left(
x,x^{\prime }\right)  \label{a1.15}
\end{equation}%
which describes essentially extremal properties of the world function $%
\sigma _{\mathrm{R}}$, i.e. the fact, that the world function $\sigma _{%
\mathrm{R}}$ generates fulfilment of condition (\ref{a1.10}).

One believes usually, that the space-time geometry is a Riemannian
(pseudo-Riemannian) geometry, and the world function $\sigma _{\mathrm{R}}$
of the real space-time geometry satisfies the equation (\ref{a1.15}). What
are reasons for such a statement?

The Riemannian space-time geometry is obtained usually as a generalization
of the proper Euclidean geometry on the case of a curved space-time
geometry. At such a generalization the original geometry is considered as a
pluralistic conception of geometry, where there are several fundamental
quantities: point, straight line, linear vector space, etc. The world
function is considered as a derivative quantity, defined by the relation (%
\ref{a1.14}). At a generalization the properties of fundamental quantities
are changed. The change of the fundamental quantities means a change of
axioms, describing properties of these quantities. The change of different
fundamental quantities is to be made by a consistent way. (New axioms of the
generalized geometry are to be compatible).

Let us compare the monistic approach and the pluralistic one in the simple
example. In the Euclidean geometry two vectors $\mathbf{P}_{0}\mathbf{P}_{1}$
and $\mathbf{Q}_{0}\mathbf{Q}_{1}$ are collinear $\left( \mathbf{P}_{0}%
\mathbf{P}_{1}\parallel \mathbf{Q}_{0}\mathbf{Q}_{1}\right) $, if and only
if the Gramm's determinant vanishes%
\begin{equation}
\mathbf{P}_{0}\mathbf{P}_{1}\parallel \mathbf{Q}_{0}\mathbf{Q}_{1}:\qquad
\left\vert 
\begin{array}{cc}
\left( \mathbf{P}_{0}\mathbf{P}_{1}.\mathbf{P}_{0}\mathbf{P}_{1}\right) & 
\left( \mathbf{P}_{0}\mathbf{P}_{1}.\mathbf{Q}_{0}\mathbf{Q}_{1}\right) \\ 
\left( \mathbf{Q}_{0}\mathbf{Q}_{1}.\mathbf{P}_{0}\mathbf{P}_{1}\right) & 
\left( \mathbf{Q}_{0}\mathbf{Q}_{1}.\mathbf{Q}_{0}\mathbf{Q}_{1}\right)%
\end{array}%
\right\vert =0  \label{b2.1}
\end{equation}%
or in the developed form%
\begin{equation}
\mathbf{P}_{0}\mathbf{P}_{1}\parallel \mathbf{Q}_{0}\mathbf{Q}_{1}:\qquad
\left( \mathbf{P}_{0}\mathbf{P}_{1}.\mathbf{Q}_{0}\mathbf{Q}_{1}\right)
^{2}=\left\vert \mathbf{P}_{0}\mathbf{P}_{1}\right\vert ^{2}\cdot \left\vert 
\mathbf{Q}_{0}\mathbf{Q}_{1}\right\vert ^{2}  \label{b2.2}
\end{equation}%
where the scalar product $\left( \mathbf{P}_{0}\mathbf{P}_{1}.\mathbf{Q}_{0}%
\mathbf{Q}_{1}\right) $ and the module $\left\vert \mathbf{P}_{0}\mathbf{P}%
_{1}\right\vert $ are expressed via world function by formulas 
\begin{equation}
\left( \mathbf{P}_{0}\mathbf{P}_{1}.\mathbf{Q}_{0}\mathbf{Q}_{1}\right)
=\sigma \left( P_{0},Q_{1}\right) +\sigma \left( P_{1},Q_{0}\right) -\sigma
\left( P_{0},Q_{0}\right) -\sigma \left( P_{1},Q_{1}\right)  \label{b2.3}
\end{equation}%
\begin{equation}
\left\vert \mathbf{P}_{0}\mathbf{P}_{1}\right\vert =\sqrt{\left( \mathbf{P}%
_{0}\mathbf{P}_{1}.\mathbf{P}_{0}\mathbf{P}_{1}\right) }=\sqrt{2\sigma
\left( P_{0},P_{1}\right) }  \label{b2.4}
\end{equation}%
The straight line $\mathcal{T}_{P_{0}P_{1},Q_{0}}$, passing through the
point $Q_{0}$ collinear to the vector $\mathbf{P}_{0}\mathbf{P}_{1}$, is
defined by the relation 
\begin{equation}
\mathcal{T}_{P_{0}P_{1},Q_{0}}=\left\{ R|\mathbf{Q}_{0}\mathbf{R\parallel P}%
_{0}\mathbf{P}_{1}\right\} =\left\{ R|\left( \mathbf{P}_{0}\mathbf{P}_{1}.%
\mathbf{Q}_{0}\mathbf{R}\right) ^{2}-\left\vert \mathbf{P}_{0}\mathbf{P}%
_{1}\right\vert ^{2}\cdot \left\vert \mathbf{Q}_{0}\mathbf{R}\right\vert
^{2}=0\right\}  \label{b2.5}
\end{equation}

If the point $Q_{0}$ coincides with the point $P_{0}$, the expression in (%
\ref{b2.5}) can be presented in the form%
\begin{equation}
\left( \mathbf{P}_{0}\mathbf{P}_{1}.\mathbf{P}_{0}\mathbf{R}\right)
^{2}-\left\vert \mathbf{P}_{0}\mathbf{P}_{1}\right\vert ^{2}\cdot \left\vert 
\mathbf{P}_{0}\mathbf{R}\right\vert ^{2}=A\left( P_{0},P_{1},R\right)
A\left( P_{0},R,P_{1}\right) B\left( P_{0},P_{1},R\right)  \label{b2.6}
\end{equation}%
where%
\begin{equation}
A\left( P_{0},P_{1},R\right) =\sqrt{\sigma \left( P_{0},R\right) }+\sqrt{%
\sigma \left( P_{1},R\right) }-\sqrt{\sigma \left( P_{1},P_{0}\right) }
\label{b2.7}
\end{equation}%
\begin{equation}
A\left( P_{0},R,P_{1}\right) =\sqrt{\sigma \left( P_{1},P_{0}\right) }+\sqrt{%
\sigma \left( P_{1},R\right) }-\sqrt{\sigma \left( P_{0},R\right) }
\label{b2.8}
\end{equation}%
\begin{equation}
B\left( P_{0},P_{1},R\right) =\sigma \left( P_{1},R\right) -\sigma \left(
P_{0},R\right) -\sigma \left( P_{1},P_{0}\right) -4\sqrt{\sigma \left(
P_{0},P_{1}\right) \sigma \left( P_{0},R\right) }  \label{b2.9}
\end{equation}%
The factor $A\left( P_{0},P_{1},R\right) $ is responsible for that part of
the straight line, which is placed between the points $P_{0}$ and $P_{1}$,
whereas the factor $A\left( P_{0},R,P_{1}\right) $ is responsible for that
part of the straight, which is placed outside the points $P_{0}$ and $P_{1}$%
. One can see from (\ref{a1.12}), that only factor $A\left(
P_{0},P_{1},R\right) $ is used in the definition of the segment $\mathcal{T}%
_{\left[ P_{0}P_{1}\right] }$.

Let us compare the monistic approach and the pluralistic one in the simple
example. In the Euclidean geometry two vectors $\mathbf{P}_{0}\mathbf{P}_{1}$
and $\mathbf{Q}_{0}\mathbf{Q}_{1}$ are the equivalent $\left( \mathbf{P}_{0}%
\mathbf{P}_{1}\text{eqv}\mathbf{Q}_{0}\mathbf{Q}_{1}\right) $, if the two
vectors are in parallel 
\begin{equation}
\mathbf{P}_{0}\mathbf{P}_{1}\uparrow \uparrow \mathbf{Q}_{0}\mathbf{Q}%
_{1}:\qquad \left( \mathbf{P}_{0}\mathbf{P}_{1}.\mathbf{Q}_{0}\mathbf{Q}%
_{1}\right) =\left\vert \mathbf{P}_{0}\mathbf{P}_{1}\right\vert \cdot
\left\vert \mathbf{Q}_{0}\mathbf{Q}_{1}\right\vert  \label{b2.10}
\end{equation}%
and their modules are equal, i.e. the following conditions are fulfilled 
\begin{equation}
\left( \mathbf{P}_{0}\mathbf{P}_{1}\text{eqv}\mathbf{Q}_{0}\mathbf{Q}%
_{1}\right) :\quad \left( \mathbf{P}_{0}\mathbf{P}_{1}.\mathbf{Q}_{0}\mathbf{%
Q}_{1}\right) =\left\vert \mathbf{P}_{0}\mathbf{P}_{1}\right\vert \cdot
\left\vert \mathbf{Q}_{0}\mathbf{Q}_{1}\right\vert \wedge \left\vert \mathbf{%
P}_{0}\mathbf{P}_{1}\right\vert =\left\vert \mathbf{Q}_{0}\mathbf{Q}%
_{1}\right\vert  \label{a1.16}
\end{equation}%
where the scalar product $\left( \mathbf{P}_{0}\mathbf{P}_{1}.\mathbf{Q}_{0}%
\mathbf{Q}_{1}\right) $ and the module $\left\vert \mathbf{P}_{0}\mathbf{P}%
_{1}\right\vert $ are expressed via world function by formulas (\ref{b2.3})
and (\ref{b2.4})

It is a well-wrought definition of two vectors equality, because it is
formulated in terms of fundamental quantity (world function), and it does
not contain a reference to means of description (coordinate system). It is
reasonable to use the definition (\ref{a1.16}) in other geometries, and in
particular, in the Riemannian geometry. In the Riemannian geometry one uses
the definition (\ref{a1.16}) only for the case, when the points $P_{0}$ and $%
Q_{0}$ coincide \ ($P_{0}=Q_{0}$). In this case the definition (\ref{a1.16})
coincides with (\ref{a1.4}). In the case, when the points $P_{0}$ and $Q_{0}$
do not coincide, the concept of the vectors $\mathbf{P}_{0}\mathbf{P}_{1}$
and $\mathbf{Q}_{0}\mathbf{Q}_{1}$ equivalence is not defined. But, why?

The answer is as follows. If in inhomogeneous geometry (in the Riemannian
geometry) the two vectors equivalence is defined by the relations (\ref%
{a1.16}), the vector equivalence appears to be multivariant, in general. It
means, that in the point $Q_{0}$ there are many vectors $\mathbf{Q}_{0}%
\mathbf{Q}_{1},\mathbf{Q}_{0}\mathbf{Q}_{1}^{\prime },\mathbf{Q}_{0}\mathbf{Q%
}_{1}^{\prime \prime },...$, which are equivalent to the vector $\mathbf{P}%
_{0}\mathbf{P}_{1}$ and are not equivalent between themselves. Multivariance
is a corollary of the fact, that to determinate the vector $\mathbf{Q}_{0}%
\mathbf{Q}_{1}$, one needs to solve two equations (\ref{a1.16}) with respect
to the point $Q_{1}$ at given points $P_{0},P_{1},Q_{0}$. In the case of the
Euclidean world function $\sigma _{\mathrm{E}}$ the solution is always
unique due to properties of the Euclidean world function $\sigma _{\mathrm{E}%
}$. However, in the case of arbitrary world function $\sigma $ there may be
many solutions, or may be no solutions at all. Multivariance of the vectors
equality leads to intransitivity of the equality relation. It means, that a
multivariant geometry is nonaxiomatizable, because in any axiomatizable
geometry the equality relation is transitive. In general, multivariance of
the equality relation leads to the fact that the straight lines appear to be
not a one-dimensional sets. It follows from the fact that the condition of
parallelism (\ref{b2.10}) is one of conditions of the vectors equivalence (%
\ref{a1.16}). The parallelism condition (\ref{b2.10}) is a special case of
the collinearity condition (\ref{b2.2}). Thus, there is a connection between
the multivariance and many-dimensionality of straight lines in multivariant
geometry.

Multivariant nonaxiomatizable geometry seems to be inadmissible for
mathematicians, who are used to consider only axiomatizable geometries. The
mathematicians consider any geometry as a logical construction, and they do
not know, how to work with a geometry, which is not a logical construction.

In the Riemannian geometry the multivariance of the equality of two vectors $%
\mathbf{P}_{0}\mathbf{P}_{1}$ and $\mathbf{Q}_{0}\mathbf{Q}_{1}$ takes
place, in general, only if $P_{0}\neq Q_{0}$. The vectors $\mathbf{P}_{0}%
\mathbf{P}_{1}$ and $\mathbf{P}_{0}\mathbf{Q}_{1}$ equality is
single-variant in the Riemannian geometry. As a result segments (\ref{a1.12})%
$\ $of geodesics (straights) in the Riemannian geometry appear to be
one-dimensional. Only straights of the form (\ref{b2.5}) appear to be
many-dimensional in the Riemannian geometry. However, such straight lines
appear practically neither in physics, nor in mathematics, and
mathematicians relate with disbelief to the statement on multivariance of a
geometry.

Using the pluralistic approach in the transition from the proper Euclidean
geometry to the Riemannian geometry, one thinks, that the one-dimensional
character of straights in the Euclidean geometry is a property of any
geometry. However, this property is a special property of the proper
Euclidean geometry, which is conditioned by the relations (\ref{a1.2}) - (%
\ref{a1.9}). In the Riemannian geometry not all properties (\ref{a1.2}) - (%
\ref{a1.9}) are fulfilled, and some of straights (geodesics) appear to
many-dimensional.

To avoid multivariance, mathematicians decided not to consider equivalence
of vectors in different points of the space. As a result in the Riemannian
geometry one considers only equivalence of vectors, having common origin.
Equivalence of vectors at different points is defined by means of special
procedure, known as a parallel transport of a vector. Procedure of the
parallel transport provides the customary single-variant equivalence of
vectors at different points, although the result depends on the path of the
parallel transport. (As a result essentially the equivalence of two vectors
at different points appears to be multivariant). Such a solution of the
vector equivalence problem corresponds to the pluralistic approach to
geometry, when in the Riemannian geometry one may change properties of
vectors in arbitrary way. As a result the Riemannian geometry appears to be
inconsistent, although this inconsistency is well masked.

At the pluralistic approach to the Riemannian geometry the world function $%
\sigma _{\mathrm{R}}$, defined by the relation (\ref{a1.12}), appears to be
derivative quantity, which depends on the shape of geodesics. Besides, in
the Riemannian geometry there may be several different geodesics, connecting
points $P_{0}$ and $P_{1}$. In this case the world function $\sigma _{%
\mathrm{R}}\left( P_{0},P_{1}\right) $, defined by the relation (\ref{a1.14}%
) appears to be many-valued. It is inadmissible in a monistic conception of
a geometry, where the world function is a fundamental quantity of a geometry.

If nevertheless we want to introduce world function $\sigma _{\mathrm{R}}$,
obtained by means of the relation (\ref{a1.14}), and to consider the
obtained world function $\sigma _{\mathrm{R}}$ as a fundamental quantity, we
need to consider only single-valued branch of $\sigma _{\mathrm{R}}$,
removing all other branches. After such a procedure one obtains the world
function as a fundamental quantity, which does not depend on choice of
geodesics. Different "Riemannian" geometries correspond to different choices
of geodesics, generating the world function.

Let a Riemannian geometry be given on the point set $\Omega $. Let us cut a
hole $\Omega _{1}$ in the set $\Omega $. At the pluralistic approach the
geometry on the set $\Omega \backslash \Omega _{1}$ changes, in general,
because the hole changes shape of several geodesics. As a result the
nonconvex point set $\Omega \backslash \Omega _{1}$, cannot be embedded
isometrically in the point set $\Omega $, because the set of geodesics,
determining the world function $\sigma _{\mathrm{R}}$ by means of (\ref%
{a1.14}), changes. It is also an inconsistency of the Riemannian geometry in
the framework of the pluralistic approach.

Using the pluralistic approach at transition from the Euclidean geometry to
the Riemannian geometry, one changes independently different fundamental
quantities. It is very difficult to change them concerted and to obtain
consistent conception of a Riemannian geometry. At the monistic approach
there is the only fundamental quantity $\sigma $. Any change of the world
function $\sigma $ generates a generalized geometry. There is no problems
with consistency of this generalized geometry, because this geometry is not
a logical construction, in general. In the obtained generalized geometry
there are no theorems and axioms, because it is a constructive geometry. All
definitions and geometrical objects of the generalized geometry are obtained
from corresponding definitions and geometrical objects of the proper
Euclidean geometry by means of a deformation. Let a statement of the
Euclidean geometry be written in terms of the Euclidean world function $%
\sigma _{\mathrm{E}}$, \textit{under condition that the special relations (%
\ref{a1.2}) - (\ref{a1.9}) of the proper Euclidean geometry are not used in
this record}. Replacing $\sigma _{\mathrm{E}}$ by the world function $\sigma 
$ of the generalized geometry in this record, one obtains the corresponding
statement of the generalized geometry. At such a deformation of the proper
Euclidean geometry one does not use the formal logic, and the problem of
inconsistency of the obtained geometry cannot be stated at all.

\textit{Important remark. }In applications to the relativity theory one uses
practically only the relation (\ref{b2.3}) for scalar product of two
vectors, which does not use the special relations (\ref{a1.2}) - (\ref{a1.9}%
).

It is of no importance, that the obtained generalized geometry is not a
logical construction. The generalized geometry is a science on shape and
disposition of geometrical objects. This circumstance is important in
applications of the geometry to physics and mechanics. It is very important,
that monistic conception of the generalized geometry does not need proofs of
numerous theorems, what is important at pluralistic approach to geometry.

Inconsistency of the Riemannian geometry, obtained in the framework of the
pluralistic approach, is a very serious balk for application of the
Riemannian geometry in the relativistic theory of gravitation. Using
monistic approach to the space-time geometry, one can generalize the general
relativity theory on the case of non-Riemannian geometry \cite{R2009,R2010}.

\section{Problems of transition from the pluralistic \newline
conception of geometry to the monistic one}

Although the monistic conception of geometry is more perfect, than the
preceding pluralistic conception, the scientific community dislikes the new
monistic conceptions as a rule. The scientific community does not
acknowledge the new monistic (or less pluralistic) conception for a long
time. Objections of the scientific community against a new monistic (or a
less pluralistic) conception have rather a social character, than a
scientific one. Indeed, reducing the number of fundamental quantities in a
new monistic conception, one is forced to transform some fundamental
quantities of the pluralistic conception into derivative quantities of the
monistic (or less pluralistic) conception. The derivative quantities of the
monistic conception look more complicated, than the aforegoing quantities of
the pluralistic conception, and members of the scientific community ask
themselves: "Why should we consider the new complicated quantities, if they
do not give anything new?"

We had such a situation in the case of transition from the axiomatic
thermodynamics to the statistical physics. Indeed, the relations of the
statistical physics and those of the kinetic theory are more complicated,
than the simple rules of work with the thermodynamic quantities. Really they
do not give anything new in those regions of physics, where the
thermodynamics works well and where most researchers work. The statistical
physics gives new results only at consideration of thermal fluctuations.
However, this circumstance was unclear for most physicists. The scientific
community as a whole was against papers by Boltzman and Gibbs, which
introduced a monistic theory of thermal phenomena, reducing thermal
phenomena to mechanical ones,.

A like situation takes place in quantum mechanics. A use of nonaxiomatizable
space-time geometry in the microcosm admits one to reduce the quantum
effects to geometrical effects. Let in microcosm the space-time geometry be
described by the world function%
\begin{equation}
\sigma _{\mathrm{d}}=\sigma _{\mathrm{M}}+d\cdot \text{sign}\sigma _{\mathrm{%
M}},\qquad d=\frac{\hbar }{2bc}=\text{const}  \label{a2.17}
\end{equation}%
where $\sigma _{\mathrm{M}}$ is the world function of the Minkowski
geometry. Here $\hbar $ is the quantum constant, $c$ is the speed of the
light, and $b$ is some universal constant. World chains (lines), describing
particle motion, are stochastic in this space-time geometry. Statistical
description of stochastic timelike world chains (lines) is equivalent to the
quantum description in terms of the Schr\"{o}dinger equation \cite{R91}. In
such a conception the quantum principles are corollaries of the space-time
geometry parameters. The number of the fundamental quantities reduces, and
the conception becomes to be less pluralistic. However, it does not obtain
any new effects in the region, where the conventional quantum mechanics
works. New effects may be obtained in the region of the elementary particles
theory \cite{R2010a}. One needs to make new investigations in the theory of
elementary particles, using the new less pluralistic conception.
Calculations, connected with these new investigations, are rather
complicated, the used formalism is new, and nobody is interested in this
less pluralistic approach until one will obtain new numerical results, which
could show, that the less pluralistic approach is valid.

The monistic conception of geometry is essential also in the megacosm. The
contemporary general relativity theory supposes, that the space-time
geometry cannot be a nonaxiomatizable geometry (a non-Riemannian geometry).
Such opinion is based on results of the contemporary geometry, where any
geometry is considered as a logical construction (but not as a science on
the dispositions of geometrical objects). For instance, the symplectic
geometry has nothing to do with properties of the space-time or of the
space. But it uses a mathematical technique, which is close to the
mathematical technique of the Euclidean geometry. Mathematicians use the
term "geometry" for the symplectic geometry, because it is a logical
construction. In other words, for mathematicians a geometry is rather a
logical construction, than a science on disposition of geometrical objects
in the space, or in the space-time.

Apparently, such a relation to geometry is a reason of abruption of
constructive (nonaxiomatizable) geometries, which are very important in
application of a geometry to physics. Mathematicians are ready to accept and
to develop inconsistent Riemannian geometry \cite{R2005}, but they are not
ready to learn monistic conception of a geometry, where everything is
copacetic with consistency. For instance, when I submitted my report
entitled "Nonaxiomatizable geometries and their application to physics" to
one of seminars of the Steklov mathematical institute, my report was
rejected. I was told only, that in the Mathematical institute there are no
researchers, who are interested in such problems. Of course, it is true.
However, the Steklov mathematical institute is a leading mathematical
institute of Russia, and its researches prefer to develop and to use the
customary inconsistent conception of a geometry, ignoring the consistent one!

Nevertheless such an approach to monistic conception of geometry is rather
reasonable in the light of the general approach to transition from a
pluralistic conception to a monistic one, although in the given case the
monistic conception of geometry is much simpler, than, for instance, the
less general Riemannian geometry.

However, there are no rules without exceptions. Such an outstanding geometer
as Grisha Perelman has evaluated the situation with monistic conception of
geometry correctly. I have not the pleasure of knowing him, and I estimate
his viewpoint, considering his behavior. After several months (end of 2005)
after publishing the paper \cite{R2005}, Perelman left the Mathematical
institute, where he worked, and said to director of the institute, that he
ceases his mathematical research. Such a decision of the prosperous geometer
can be explained only by a strong disappointment in his activity. During all
his life Perelman developed topological direction in geometry, where
topological quantities are the fundamental quantities of geometry. When he
has understood, that the topological quantities are derivative quantities,
which are determined by the world function, and the Riemannian geometry,
which has been used in his investigations, is inconsistent, he was shocked.
Apparently, this shock determines his decision on termination of his
mathematical investigations.

As concerns to his refuse from the Fields medal and from award of the Clay
Mathematics Institute, it is important only from the viewpoint of mercantile
grassroots, who seek material values instead of human. From viewpoint of
Perelman any money are nothingness with respect to problems of mathematics.
Besides, the problem, solved by Perelman, is a Millennium Prize Problems
only now. After several years this problem will transform to usual problem
of topology. Perelman understood this very well. He refused from the award,
referring to incompetence of people, adjudging a prize to him. Indeed, one
can understand behavior of Perelman, only if one does know about the
Riemannian geometry inconsistency. But nobody pay attention onto this
inconsistency, even if one knows about signs of this inconsistency.

I have not discussed with Perelman his relation to the Riemannian geometry.
My interpretation of the Perelman's behaviour is only a hypothesis. However,
it is a very verisimilar hypothesis. I do not now other reasonable
hypothesis, which could explain the Perelman's behaviour. I have presented
this hypothesis, to convince those researchers, who believe rather in
authorities, than in logical arguments. The Riemannian geometry in its
conventional presentation is inconsistent, even if my interpretation of the
Perelman's behaviour is wrong.

Thus, the relation of researchers to the monistic conception of geometry
develops in its natural way. I think, that those researchers, who have
understood advantages of the monistic conception of geometry and will use it
in their investigations, will succeed in progress of physics.

\end{document}